\documentclass[10pt, reqno]{amsart}

\usepackage{amssymb}
\usepackage{bm}
\usepackage{graphicx}
\usepackage{xcolor}
\usepackage{hyperref}
\usepackage{mathdots}

\numberwithin{equation}{section}

\hypersetup{
	colorlinks=true,
	linkcolor=blue,
	citecolor=green
	}

\theoremstyle{plain}
\newtheorem{thm}{Theorem}[section]

\newtheorem{conj}[thm]{Conjecture}

\theoremstyle{definition}

\newtheorem{exa}[thm]{Example}

\theoremstyle{remark}
\newtheorem{rmk}[thm]{Remark}

\newcommand{\R}{\mathbb{R}}
\newcommand{\Z}{\mathbb{Z}}
\newcommand{\N}{\mathbb{N}}
\newcommand{\Mcal}{\mathcal{M}}
\newcommand{\Ical}{\mathcal{I}}
\newcommand{\Nfra}{\mathfrak{N}}
\newcommand{\ga}{\bm{\gamma}}
\newcommand{\ee}{\mathrm{e}}
\newcommand{\ii}{\mathrm{i}}
\newcommand{\dd}{\mathrm{d}}

\newcommand{\OO}{\mathrm{O}}

\begin{document}
	
	\title[counting integral points near curves]{counting integral points near curves: an exposition}
	
	\author{Tianliang Wu}
	\address{School of Mathematical Sciences, University of Science and Technology of China, Hefei, 230026, P.R. China}
	\email{wutianliang@mail.ustc.edu.cn}
	
	\subjclass[2020]{Primary 11D75; Secondary 11P21}
	
	\keywords{Integral points near curves, lattice point counting}
	
	\begin{abstract}
		We give a detailed account of the recent work of Hickman and Srivastava~\cite{hickman2025counting} on counting integral points in a $ \delta $-neighbourhood of a dilate $ q \Mcal $ of a non-degenerate curve $ \Mcal \subset \R^{n} $, $n\ge 3$. Our aim is expository: we simplify the proof, provide additional details, refine the statement of the theorem and provide calculations from a more fundamental perspective. We also formulate a new conjecture based on recent work of Chen, Seeger, Srivastava and Technau \cite[Lemma~2.1]{chen2026sharp}.
	\end{abstract}
	
	\maketitle
	
	\section{Introduction}
	
	Let $ \Mcal $ be a $d$-dimensional compact submanifold of $ \R^n $, $ 1 \leq d < n $, and let $ A_\Mcal(q,\delta) $
	denote the cardinality of the set
	\[
		\{ \bm{b} \in \Z^n: \operatorname{dist} (\bm{b} / q, \Mcal) \leq \delta / q \} = 
		\{ \bm{b} \in \Z^n: \operatorname{dist} (\bm{b}, q\Mcal) \leq \delta \},
	\]
	where $ q \geq 1 $, $ \delta \geq 0 $, and $ \operatorname{dist} $ denotes the distance with respect to the Euclidean norm on $ \R^n $.
	
	Throughout the paper, $ A \ll_L B $ or $ A = \OO_L(B) $ (where $ L $ is a list of variables and could be empty, and $ A $ and $ B $ are expressions) means there exists a constant $ C_L > 0 $ depending only on $ \Mcal $, $ n $ and the variables in $ L $ such that $ A \leq C_L B $. Also, $ A \gg_L B $ means $ B \ll_L A $, and $ A \asymp_L B $ means both $ A \ll_L B $ and $ A \gg_L B $.
	
	\subsection{Background}
	Over the years, researchers have obtained upper and lower bounds for the counting function $ A_\Mcal(q,\delta) $ by imposing various conditions on $ \Mcal $.
	
	If we fix $ \delta = 0 $, the function $ A_\Mcal(q,0) $ counts the number of integral points on the dilated manifold $ q\Mcal $. Results on estimating $ A_\Mcal(q,0) $ have primarily focused on the case where $ \Mcal $ is a planar curve, i.e., $ d=1 $, $ n=2 $. In 1974, Swinnerton-Dyer \cite{swinnerton1974number} showed that if $ \Mcal $ is a $ \mathrm{C}^3 $ strictly convex compact planar curve, then
	\[
		A_\Mcal(q,0) \ll_\epsilon q^{\frac{3}{5}+\epsilon}.
	\]
	Several years later, Schmidt \cite{schmidt1985integer} gave a uniform (in $ \Mcal $) version of Swinnerton-Dyer's theorem and generalized it to hypersurfaces (i.e., $ n-d=1 $) satisfying certain curvature conditions. Soon after, Bombieri and Pila \cite{bombieri1989number} developed an ingenious determinant method and proved the essentially optimal bound
	\[
		A_\Mcal(q,0) \ll_\epsilon q^{\frac{1}{2}+\epsilon},
	\]
	where $ \Mcal $ is a smooth strictly convex compact planar curve.
	
	For $ \delta > 0 $, the function $ A_\Mcal(q,\delta) $ counts the number of integral points in the $ \delta $-neighbourhood of the dilated manifold $ q\Mcal $. In general, the problem of estimating $ A_\Mcal(q,\delta) $ becomes more difficult as the codimension $ n - d $ of the submanifold $ \Mcal $ increases. For the case $ d=1 $, $ n=2 $, various estimates and asymptotic formulas are available; see, for example, \cite{huxley1996area, trifonov2002lattice, huang2020diophantine}. For the case $ n - d = 1, n \geq 3 $, Lettington \cite{lettington2009integer, lettington2010integer} proved essentially best possible bounds for convex hypersurfaces under suitable curvature conditions. For the case $ d=1 $, $ n=3 $, Huang \cite{huang2019integral} made the first nontrivial progress on estimating $ A_\Mcal(q,\delta) $. Specifically, he proved the following theorem:
	
	\begin{thm}[\cite{huang2019integral}, Theorem 1] \label{thm:huang}
		Let $ \Mcal $ be a $ \mathrm{C}^3 $ curve in $ \R^3 $ with nonvanishing torsion. Then for all $ \delta \in [0,1/2) $ and $ q \geq 2 $, we have
		\[
			A_\Mcal(q,\delta) \ll \delta^2 q + q^\frac{3}{5} (\log{q})^\frac{4}{5}.
		\]
		Furthermore, for all $ q \gg 1 $ and $ \delta \gg q^{-\frac{1}{5}} (\log{q})^\frac{2}{5} $, we have
		\[
			A_\Mcal(q,\delta) \gg \delta^2 q.
		\]
	\end{thm}
	
	\begin{rmk}
		The sentence ``\,for all $ q \gg 1 $ and $ \delta \gg q^{-\frac{1}{5}} (\log{q})^\frac{2}{5} $\,'' should be understood as ``\,there exist positive constants $ Q_0 $ and $ C_0 $ depending only on $ \Mcal $ such that for all $ q > Q_0 $ and $ \delta > C_0 q^{-\frac{1}{5}}  (\log{q})^\frac{2}{5} $\,''.
	\end{rmk}
	
	A few years after Huang's result, Hickman and Srivastava \cite{hickman2025counting} obtained the first nontrivial upper and lower bounds for the case $ d = 1 $, $ n \geq 4 $. We discuss their theorem in detail below. There are also results for special classes of manifolds. For example, Huang and Liu \cite{huang2019simultaneous} derived optimal results for affine subspaces of $ \R^n $ satisfying certain Diophantine type conditions.
	
	Additionally, one may sum $ A_\Mcal(q,\delta) $ over $ q $ and estimate $ \sum_{q=1}^Q A_\Mcal(q,\delta) $ for $ Q \in \N $ and $ \delta > 0 $. This is called the problem of counting rational points near the manifold $ \Mcal $ and is closely related to questions in Diophantine approximation and the dimension growth problem.
	Estimating the number of rational points near manifolds is a very active area of research, and we refer the interested reader to \cite{beresnevich2012rational, huang2020density, beresnevich2023khintchine, chen2025rational, srivastava2025counting, chen2026sharp} and the references therein.
	
	\subsection{Example and conjecture}
	For the rest of the paper we assume $ \Mcal $ to be a compact $ \mathrm{C}^\infty $ curve in $ \R^n $, where $ n \geq 3 $. We say that $ \Mcal $ is \emph{non-degenerate} if for any local parametrization $ \{ \ga(t) : t \in U \} $ of $ \Mcal $, we have
	\[
		\det \left[\ga^{(1)}(t),\ga^{(2)}(t),\dots,\ga^{(n)}(t)\right] \neq 0 \quad \text{for any } t \in U.
	\]
	
	The non-degeneracy condition guarantees that the curve is sufficiently curved so that it is not locally contained in any proper linear subspace of $ \R^n $. We will always assume $ \Mcal $ to be non-degenerate so that we can obtain nontrivial estimates on $ A_\Mcal(q,\delta) $.
	
	\begin{exa}
		Let $ \Mcal = \{ (\sqrt{3} \cos{t}, \sqrt{3} \sin{t}, t) : t \in [0,2\pi] \} \subset \R^3 $, which is obviously non-degenerate. Since there are no integral points on the planar curve $ x^2 + y^2 = 3 q^2 $ for any $ q \in \N $, there are no integral points on $ q\Mcal $. Then by the compactness of $ \Mcal $ we know that $ A_\Mcal(q,\delta) = 0 $ for all $ q \in \N $ and small $ \delta $ (the threshold value depends on $ q $).
	\end{exa}
	
	The preceding example shows that the counting function may vanish for sufficiently small $ \delta $. To get a nontrivial lower bound, one might have to assume that $ \delta $ is greater than a threshold value which depends on $ q $.
	
	The following discussion is heuristic and is intended to motivate the form of the lower and upper bounds. Intuitively, if we denote the volume of the $ \delta $-neighbourhood of $ q\Mcal $ by $ V_\Mcal(q,\delta) $, then $ V_\Mcal(q,\delta) \asymp \delta^{n-1} q $. We define the error term
	\[
		E_\Mcal(q,\delta) = | A_\Mcal(q,\delta) - V_\Mcal(q,\delta) |,
	\]
	and two sets
	\[
		I_1(q) = \left\lbrace  \delta \in \left[ 0,\frac12 \right) : E_\Mcal(q,\delta) \leq \frac12 V_\Mcal(q,\delta) \right\rbrace, \, I_2(q) = \left[ 0,\frac12 \right) \setminus I_1(q).
	\]
	Set $ \delta^*(q) = \sup I_2(q) $. In practice, we usually have
	\[
		E_\Mcal(q,\delta^*(q)) \asymp V_\Mcal(q,\delta^*(q)) \asymp \delta^*(q)^{n-1} q.
	\]
	Now if $ \delta \in I_1(q) $, then $ A_\Mcal(q,\delta) \ll \delta^{n-1} q $. If $ \delta \in I_2(q) $, then by the monotonicity of $ A_\Mcal(q,\delta) $ in $ \delta $, we have
	\[
		A_\Mcal(q,\delta) \leq A_\Mcal(q,\delta^*(q)) \leq V_\Mcal(q,\delta^*(q)) + E_\Mcal(q,\delta^*(q))
		\ll \delta^*(q)^{n-1} q.
	\]
	So we get the upper bound
	\[
		A_\Mcal(q,\delta) \ll \delta^{n-1} q + \delta^*(q)^{n-1} q
		\quad \text{for } q \geq 1 \text{ and } \delta \in \left[ 0,\frac12 \right).
	\]
	Also, it is easy to get the lower bound
	\[
		A_\Mcal(q,\delta) \gg \delta^{n-1} q
		\quad \text{for } q \geq 1 \text{ and } \frac12 > \delta > \delta^*(q).
	\]

	In the case of the next example, we have $ \delta^*(q) \gg q^{-\frac1n} $ for infinitely many $ q $. So it is conjectured that
	\[
		A_\Mcal(q,\delta) \ll \delta^{n-1} q + q^\frac1n \quad \text{for } q \geq 1 \text{ and } \delta \in \left[ 0, \frac12 \right),
	\]
	and
	\[
		A_\Mcal(q,\delta) \gg \delta^{n-1} q \quad \text{for } q \geq 1 \text{ and } \frac12 > \delta \gg q^{-\frac1n}.
	\]
	
	\begin{exa}
		Let $ \Mcal = \{ (t,t^2,\dots,t^n) : t \in [0,1] \} \subset \R^n $, which is obviously non-degenerate. Suppose that $ q = k^n $ for some $ k \in \N $. Then there are exactly $ k + 1 $ integral points on $ q\Mcal $ corresponding to the parameter values $ t = 0,\frac{1}{k},\frac{2}{k},\dots,1 $. Therefore,
		\[
			A_\Mcal(q,\delta) \geq A_\Mcal(q,0) = q^\frac{1}{n} + 1 \asymp q^\frac{1}{n}
		\]
		for all $ q = k^n $, $ k \in \N $, $ \delta > 0 $.
	\end{exa}
	
	However, recently Chen, Seeger, Srivastava and Technau \cite[Lemma~2.1]{chen2026sharp} proved that for $ \epsilon > 0 $ and $ \Mcal $ as in the last example, we have
	\[
		\sum_{q=1}^Q A_\Mcal(q,\delta) \gg \delta^\frac12 Q^\frac32 \quad \text{for } \frac12 > \delta > Q^{\epsilon - 1} \text{ and } Q \gg_\epsilon 1.
	\]
	So the above conjecture is incorrect for $ n \geq 4 $. The best we can hope for (if we only consider the bounds of the form ``\,$ \delta^{n-1} q + q^{\theta(n)} $\,'') is the following conjecture:
	\begin{conj} \label{conj}
		Let $ \Mcal $ be a compact $ \mathrm{C}^\infty $ non-degenerate curve in $ \R^n $, $ n \geq 3 $. Then
		\[
			A_\Mcal(q,\delta) \ll \delta^{n-1} q + q^\frac{n-2}{2n-3} \quad \textup{for } q \geq 1 \textup{ and } \delta \in \left[ 0, \frac12 \right),
		\]
		and
		\[
			A_\Mcal(q,\delta) \gg \delta^{n-1} q \quad \textup{for } q \geq 1 \text{ and } \frac12 > \delta \gg q^{-\frac{1}{2n-3}}.
		\]
	\end{conj}
	
	In fact, if $ A_\Mcal(q,\delta) \ll \delta^{n-1} q + q^{\theta} $, then $ \sum_{q=1}^Q A_\Mcal(q,\delta) \ll \delta^{n-1} Q^2 + Q^{1 + \theta} $. So we should have $ \delta^\frac12 Q^\frac32 \ll \delta^{n-1} Q^2 + Q^{1 + \theta} $, which is equivalent to
	\[
		\delta^{n-\frac32} Q^\frac12 + \delta^{-\frac12} Q^{-\frac12 + \theta} \gg 1.
	\]
	Then we have $ \theta \geq (n-2) / (2n-3) $.
	
	\subsection{The theorem of Hickman–Srivastava}
	Hickman and Srivastava obtained the following results for $ n \geq 3 $. For $ n = 3 $, this theorem recovers the bounds of Theorem~\ref{thm:huang} up to $ \varepsilon $ losses.
	
	\begin{thm}[Hickman-Srivastava \cite{hickman2025counting}, with some refinements] \label{thm:hickman&srivastava}
		For $ n \geq 3 $, let $ \Mcal $ be a compact $ \mathrm{C}^\infty $ non-degenerate curve in $ \R^n $ and
		\[
			\Theta(n) = \begin{cases} 
				\frac{n^2+4}{n(n+4)} & \textup{for even } n, \\ 
				\frac{n^2+3}{n(n+4)-1} & \textup{for odd } n. 
			\end{cases}
		\]
		Let $ \nu > 0 $ and $ \Delta > 0 $. Then for all $ \delta \in [0,\Delta) $ and $ q \geq 1 $, we have the upper bound
		\[
			A_\Mcal(q,\delta) \ll_{\nu,\Delta} \delta^{n-1} q + q^{\Theta(n) + \nu}.
		\]
		Further, for all $ q \gg_\nu 1 $ and $ \delta \gg q^{(\Theta(n)-1)/(n-1) + \nu} $, we also have the lower bound
		\[
			A_\Mcal(q,\delta) \gg \delta^{n-1} q.
		\]
	\end{thm}
	
	\begin{rmk}
		The sentence ``\,for all $ q \gg_\nu 1 $ and $ \delta \gg q^{(\Theta(n)-1)/(n-1) + \nu} $\,'' should be understood as ``\,there exist a positive constant $ Q_0 $ depending only on $ \Mcal $, $ n $ and $ \nu $, and a positive constant $ C_0 $ depending only on $ \Mcal $ and $ n $, such that for all $ q > Q_0 $ and $ \delta > C_0 q^{(\Theta(n)-1)/(n-1) + \nu} $\,''.
	\end{rmk}
	
	\begin{rmk} \label{remark}
		The ``\,$ \mathrm{C}^\infty $\,'' regularity assumption can be relaxed to ``\,$ \mathrm{C}^{2n} $\,'' regularity, since we have used derivatives of order at most $ 2n $ in the subsequent proof.
	\end{rmk}
	
	\subsection{What is new in this note}
	\begin{enumerate}
		\item \textbf{Simplification.}
		We avoid splitting the sum into near and far parts. When ``\,reversing\,'' the Frenet frame, we start with $ \bm{e}_n $ rather than $ \bm{e}_n / e_{n,k} $.
		\item \textbf{Details supplied.}
		We provide a detailed explanation (in Section~\ref{section:reductions and smoothing}) of why the curve $ \Mcal $ can be reduced to the standard form (\ref{def:standard form}).
		\item \textbf{New conjecture.}
		We formulate Conjecture~\ref{conj} based on recent work of Chen, Seeger, Srivastava and Technau \cite[Lemma~2.1]{chen2026sharp}.
		\item \textbf{Different calculation.}
		In Section~\ref{section:final calculation}, we provide calculations from a more fundamental perspective, from which one may understand that the estimates in \cite{hickman2025counting} are actually sharp within the framework of this method. To get better results, one needs to apply different methods.
		\item \textbf{Refinement of the theorem.} As shown in Theorem~\ref{thm:hickman&srivastava}, the admissible range of $ \delta $ is extended to $ [0,\Delta) $, and in the lower bound $ A_\Mcal(q,\delta) \gg \delta^{n-1} q $, the implicit constant is independent of $ \nu $. Also in Remark~\ref{remark} we point out that the ``\,$ \mathrm{C}^\infty $\,'' regularity assumption can be relaxed to ``\,$ \mathrm{C}^{2n} $\,'' regularity.
	\end{enumerate}
	
	Everything else in this note is due to Hickman and Srivastava; we have tried to attribute carefully, and any inaccuracy in the exposition is entirely our own.
	
	\section{Reduction and smoothing} \label{section:reductions and smoothing}
	
	\subsection{Reduction to standard form}
	Recall that $ \Mcal $ is a compact non-degenerate smooth curve in $ \R^n $. By the inverse function theorem, any point on $ \Mcal $ has a neighbourhood which can be parametrized in one of the following $ n $ forms:
	\[
		\{ (f_{1}(t),\dots,f_{i-1}(t),t,f_{i+1}(t),\dots,f_{n}(t)) : t \in [a,b] \}, \quad i=1,2,\dots,n.
	\]
	We shrink each of these pieces of $ \Mcal $, by restricting the parameter interval from $ [a,b] $ to $ \frac12 [a,b] = [a+\frac{b-a}{4},b-\frac{b-a}{4}] $.
	
	Since $ \Mcal $ is compact, it can be covered by finitely many such smaller pieces. Therefore, counting integral points near the dilation of $ \Mcal $ reduces to counting integral points near the dilation of each of those smaller pieces. Also, by relabeling coordinates, which does not affect the lattice point counting, we may assume each of those smaller pieces has the form
	\[
		\{ (t,f_{2}(t),f_{3}(t),\dots,f_{n}(t)) : t \in \frac12 [a,b] \}.
	\]
	
	For simplicity, we may assume that $ \Mcal $ has the above form and fix the parameter interval, that is, we assume $ \Mcal $ has the \emph{standard form}
	\begin{equation}
		\Mcal = \{ \ga(t) \stackrel{\text{def}}{=} (t,\bm{f}(t)) : t \in I \stackrel{\text{def}}{=} [-1,1] \}, \label{def:standard form}
	\end{equation}
	where $ \bm{f} = (f_{2},\dots,f_{n}) $ is smooth and defined on $ 2I=[-2,2] $, and $ \ga^{(1)}(t),\dots,\ga^{(n)}(t) $ are linearly independent for any $ t \in 2I $.
	
	In the rest of this paper we will prove Theorem~\ref{thm:hickman&srivastava} only for the case in which $ \Mcal $ has the standard form. As one can check, with some minor modifications, our proof works equally well if we replace $ I $ by any compact interval $ [a,b] $ with $ b>a $. Thus, by the above argument, we will have actually proved Theorem~\ref{thm:hickman&srivastava} for general $ \Mcal $.
	
	\subsection{Another counting function}
	For the rest of the paper we assume that $ \Mcal $ has the standard form (\ref{def:standard form}). Define another version of the counting function
	\[
		N_\Mcal(q,\delta) = \# \left\lbrace  a \in \Z : \frac{a}{q} \in I, \; \left\| q f_i \left( \frac{a}{q} \right) \right\| \leq \delta, \; i=2,\dots,n \right\rbrace,
	\]
	where $ \| \cdot \| = \operatorname{dist}(\, \cdot \, ,\Z) $, $ q \geq 1 $ and $ \delta \geq 0 $.
	
	Note that $ N_\Mcal(q,\delta) = 2 \lfloor q \rfloor + 1 $ when $ \delta \geq 1/2 $. Moreover, it is not hard to prove that for all $ q \geq 1 $ and $ \delta \in (0,\Delta) $, we have
	\begin{equation}
		N_\Mcal(q,\delta/n) \leq A_\Mcal(q,\delta) \ll_\Delta 1 + N_\Mcal (q, C_\Mcal \delta), \label{ineq:A}
	\end{equation}
	where $ C_\Mcal = \max_{2 \leq i \leq n, \, t \in [-1,1]} |f_i^{'}(t)| + 1 $. Then we reduce Theorem~\ref{thm:hickman&srivastava} to the following theorem:
	
	\begin{thm} \label{thm:reduced}
		For $ n \geq 3 $, let $ \Mcal $ be of the standard form (\ref{def:standard form}). Let $ \nu > 0 $. Then for all $ \delta \in (0,1/2) $ and $ q \geq 1 $, we have the upper bound
		\[
			N_\Mcal(q,\delta) \ll_\nu \delta^{n-1} q + q^{\Theta(n) + \nu}.
		\]
		Further, for all $ q \gg_\nu 1 $ and $ 1/2 > \delta \gg q^{(\Theta(n)-1)/(n-1) + \nu} $, we also have the lower bound
		\[
			N_\Mcal(q,\delta) \gg \delta^{n-1} q.
		\]
	\end{thm}
	
	\subsection{Smooth counting function}
	For the rest of the paper, we need only to prove Theorem~\ref{thm:reduced}. To estimate the counting function $ N_\Mcal $, we introduce a smoothed version.
	
	Let $ \omega^+ \in \mathrm{C}_\mathrm{c}^\infty(\R) $ take values in $ [0,1] $ such that $ \operatorname{supp}{\omega^+} \subset [-2,2] $ and $ \omega^+(t) = 1 $ for $ t \in [-1,1] $. Define $ \omega^- \in \mathrm{C}_\mathrm{c}^\infty(\R) $ by $ \omega^-(t) = \omega^+(2t) $.
	
	Similarly, let $ \eta^+ \in \mathrm{C}_\mathrm{c}^\infty(\R^{n-1}) $ take values in $ [0,1] $ such that $ \operatorname{supp}{\eta^+} \subset [-2,2]^{n-1} $ and $ \eta^+(t) = 1 $ for $ t \in [-1,1]^{n-1} $. Define $ \eta^- \in \mathrm{C}_\mathrm{c}^\infty(\R^{n-1}) $ by $ \eta^-(t) = \eta^+(2t) $.
	
	For $ \omega \in \{ \omega^+,\omega^- \} $, $ \eta \in \{ \eta^+,\eta^- \} $, $ q \geq 1 $ and $ \delta \in (0,1/2) $, define the smooth counting function
	\[
		\Nfra_{\omega,\eta,\Mcal}(q,\delta) = \sum_{a \in \Z} \sum_{\bm{b} \in \Z^{n-1}}
		\omega \left( \frac{a}{q} \right)
		\eta \left( \frac{1}{\delta} \left(  q \bm{f} \left( \frac{a}{q} \right) - \bm{b} \right)  \right).
	\]
	It is easy to see that for all $ q \geq 1 $ and $ \delta \in (0,1/2) $, we have
	\begin{equation}
		\Nfra_{\omega^-,\eta^-,\Mcal}(q,\delta)
		\leq N_\Mcal(q,\delta)
		\leq \Nfra_{\omega^+,\eta^+,\Mcal}(q,\delta). \label{ineq:N}
	\end{equation}
	Thus, what remains is to estimate $ \Nfra_{\omega^-,\eta^-,\Mcal} $ and $ \Nfra_{\omega^+,\eta^+,\Mcal} $. We will abbreviate $ \Nfra_{\omega,\eta,\Mcal} $ as $ \Nfra_\Mcal $ when there is no ambiguity.
	
	In addition, by choosing $ \eta^+ $ appropriately, we may assume that for fixed $ q $, the function $ \Nfra_\Mcal(q,\delta) $ is non-decreasing in $ \delta $.
	
	\subsection{Poisson summation}
	We are now able to use Fourier tools to estimate the smooth counting function $ \Nfra_\Mcal $. We give only a brief outline here; see \cite[Section~2]{hickman2025counting} for details.
	
	The Poisson summation formula gives
	\[
		\sum_{\bm{b} \in \Z^{n-1}} \eta \left( \delta^{-1} (\bm{x} - \bm{b}) \right)
		= \delta^{n-1} \sum_{\bm{j}_0 \in \Z^{n-1}}
		\widehat{\eta}(\delta \bm{j}_0)
		\ee^{2 \pi \ii \bm{x} \cdot \bm{j}_0},
		\quad \bm{x} \in \R^{n-1},
	\]
	where $ \widehat{\eta}(\bm{\xi}) = \int_{\R^{n-1}} \eta(\bm{x}) \ee^{-2\pi\ii\bm{x}\cdot\bm{\xi}} \, \dd \bm{x} $ is the Fourier transform of $ \eta $. Then we can represent $ \Nfra_\Mcal $ in terms of an exponential sum
	\[
		\Nfra_\Mcal(q,\delta) = \delta^{n-1} \sum_{a \in \Z} \sum_{\bm{j}_0 \in \Z^{n-1}}
		\widehat{\eta}(\delta \bm{j}_0)
		\ee^{2 \pi \ii q \bm{f}(a/q) \cdot \bm{j}_0} \omega (a/q).
	\]
	Using the Schwartz decay of $ \widehat{\eta} $, for any $ \varepsilon > 0 $, we have
	\begin{align*}
		\sum_{\bm{j}_0 \in \Z^{n-1}}
		\widehat{\eta}(\delta \bm{j}_0) \ee^{2 \pi \ii q \bm{f}(a/q) \cdot \bm{j}_0}
		&= \sum_{\| \bm{j}_0 \|_\infty \leq q^\varepsilon \delta^{-1}}
		\widehat{\eta}(\delta \bm{j}_0) \ee^{2 \pi \ii q \bm{f}(a/q) \cdot \bm{j}_0} \\
		&\quad + \OO_\varepsilon \left( \delta^{-(n-1)}q^{-1} \right) .
	\end{align*}
	Combining the previous two equations and isolating the term $ \bm{j}_0 = \bm{0} $, we get
	\begin{align*}
		\Nfra_\Mcal(q,\delta)
		&= c_0(q) \delta^{n-1} q \\
		&\quad + \delta^{n-1}
		\sum_{1 \leq \| \bm{j}_0 \|_\infty \leq q^\varepsilon \delta^{-1}} \widehat{\eta}(\delta \bm{j}_0)
		\sum_{a \in \Z} \ee^{2 \pi \ii q \bm{f}(a/q) \cdot \bm{j}_0} \omega (a/q) \\
		&\quad + \OO_\varepsilon(1),
	\end{align*}
	where $ c_0(q) = \widehat{\eta}(0) q^{-1} \sum_{a \in \Z} \omega (a/q) \asymp 1 $.
	
	Using Poisson summation again, along with non-stationary phase, for $ \| \bm{j}_0 \|_\infty \leq q^\varepsilon \delta^{-1} $, we have (see \cite[Lemma~2.1]{hickman2025counting} for details)
	\begin{align*}
		\sum_{a \in \Z} \ee^{2 \pi \ii q \bm{f}(a/q) \cdot \bm{j}_0} \omega (a/q)
		&= q \sum_{k \in \Z}
		\int_\R \ee^{2 \pi \ii q (\bm{f}(t) \cdot \bm{j}_0 + t k)} \omega(t) \, \dd t \\
		&= q \sum_{|k| \leq M_{\ga} q^\varepsilon \delta^{-1}}
		\int_\R \ee^{2 \pi \ii q (\bm{f}(t) \cdot \bm{j}_0 + t k)} \omega(t) \, \dd t \\
		&\quad + \OO_\varepsilon(q^{-(n-1)\varepsilon}),
	\end{align*}
	where $ M_{\ga} = 2 \sup \{ \| \partial_t \bm{f}(t) \|_1 : t \in [-2,2] \} $. Now we introduce the notation
	\[
		\Ical_\Mcal(\bm{\xi}) = \int_\R \ee^{2 \pi \ii \ga(t) \cdot \bm{\xi}} \omega(t) \, \dd t,
		\quad \bm{\xi} \in \R^n.
	\]
	Write $ \bm{j} = (k,\bm{j}_0) $, and recall that $ \ga(t) = (t,\bm{f}(t)) $, so
	\begin{align}
		\Nfra_\Mcal(q,\delta)
		&= c_0(q) \delta^{n-1} q \notag \\
		&\quad + \delta^{n-1} q \sum_{\substack{k \in \Z \\ |k| \leq M_{\ga} q^\varepsilon \delta^{-1}}} \sum_{\substack{\bm{j}_0 \in \Z^{n-1} \\ 1 \leq \| \bm{j}_0 \|_\infty \leq q^\varepsilon \delta^{-1}}} \widehat{\eta}(\delta \bm{j}_0) \, \Ical_\Mcal(q \bm{j}) \notag \\
		&\quad + \OO_\varepsilon(1). \label{eq:Ncal}
	\end{align}
	
	\section{The H-functional} \label{section:the H-functional}
	Let $ \ga $ be as in (\ref{def:standard form}), and define the \emph{H-functional}
	\[
		H(\bm{\xi}) = \inf_{t \in [-2,2]} \max_{1 \leq r \leq n} |\ga^{(r)}(t) \cdot \bm{\xi}|^\frac{1}{r}, \quad \bm{\xi} \in \R^n.
	\]
	Integration by parts and an Arkhipov–Chubarikov–Karatsuba-type estimate (see \cite[Section~3.1]{hickman2025counting} for details) yield that for all $ |\bm{\xi}| \geq 1 $, we have
	\[
		|\Ical_\Mcal(\bm{\xi})| \ll H(\bm{\xi})^{-1}.
	\]
	Then by (\ref{eq:Ncal}) and the triangle inequality, we get
	\begin{equation}
		|\Nfra_\Mcal(q,\delta) - c_0(q) \delta^{n-1} q| \ll \delta^{n-1} q
		\sum_{\substack{\bm{j} \in \Z^n \\ 1 \leq | \bm{j} | < C q^\varepsilon \delta^{-1}}}
		H(q \bm{j})^{-1} + \OO_\varepsilon(1), \label{ineq:Nfra}
	\end{equation}
	where $ C $ is a sufficiently large constant depending only on $ \ga $ (or equivalently, $ \Mcal $) and $ n $.
	
	\subsection{Level set decomposition} By (\ref{ineq:Nfra}), all we need to do is to estimate 
	\[
		\varSigma_\Mcal(q,\delta) \stackrel{\text{def}}{=} \delta^{n-1} q
		\sum_{\substack{\bm{j} \in \Z^n \\ 1 \leq | \bm{j} | < C q^\varepsilon \delta^{-1}}}
		H(q \bm{j})^{-1}.
	\]
	We use a dyadic level set decomposition.
	
	Using the non-degeneracy condition, it is elementary to prove that
	\[
		|\bm{\xi}|^\frac{1}{n} \ll H(\bm{\xi}) \ll |\bm{\xi}|, \quad |\bm{\xi}| \geq 1.
	\]
	Enlarging the constant $ C $ if necessary, we may assume that
	\[
		2 C^{-1} |\bm{\xi}|^\frac{1}{n} < H(\bm{\xi}) < \frac12 C |\bm{\xi}|, \quad |\bm{\xi}| \geq 1.
	\]
	Fix $ q \geq 1 $ and define the level set
	\[
		S(R,\lambda) = \{ \bm{j} \in \Z^n : R \leq |q \bm{j}| < 2R, \, \lambda \leq H(q \bm{j}) < 2 \lambda \},
	\]
	where $ R \geq q $ and $ \lambda \geq C^{-1} $.
	Thus we can decompose the sum as follows:
	\begin{align}
		\varSigma_\Mcal(q,\delta)
		&\leq \delta^{n-1} q
		\sum_{\substack{1 \leq q^{-1} R < C q^\varepsilon \delta^{-1} \\ q^{-1} R \; \text{dyadic}}}
		\sum_{\substack{\bm{j} \in \Z^n \\ R \leq |q \bm{j}| < 2R}} H(q \bm{j})^{-1} \notag \\
		&\leq \delta^{n-1} q
		\sum_{\substack{1 \leq q^{-1} R < C q^\varepsilon \delta^{-1} \\ q^{-1} R \; \text{dyadic}}}
		\sum_{\substack{C^{-1} R^\frac{1}{n} < \lambda < C R \\ \lambda \; \text{dyadic}}}
		\# S(R,\lambda) \cdot \lambda^{-1}, \label{ineq:Sigma}
	\end{align}
	where ``\,$ q^{-1} R \; \text{dyadic} $\,'' means $ q^{-1} R $ takes values in $ \{ 2^m : m \in \Z \} $, and similarly for ``\,$ \lambda \; \text{dyadic} $\,''.
	
	\subsection{Geometry of the level set}
	To estimate the cardinality of the level set $ S(R,\lambda) $, we need a detailed geometric description of it. For this purpose, we exploit the Frenet frame of the non-degenerate curve $ \Mcal $. For simplicity, we consider some larger set
	\begin{multline*}
		S(R,\lambda) \subset S_1(R,\lambda) \stackrel{\text{def}}{=} \{ \bm{j} \in \Z^n : |q \bm{j}| < 2R, \, H(q \bm{j}) < 2 \lambda \} \\
		= \bigcup_{t \in [-2,2]} \{ \bm{j} \in \Z^n : |q \bm{j}| < 2R, \, |\bm{j} \cdot \ga^{(r)}(t)| < q^{-1} (2 \lambda)^r, \, r=1,\dots,n \}.
	\end{multline*}

	Applying the Gram-Schmidt process to the vectors $ \ga^{(1)},\dots,\ga^{(n)} $, we get the Frenet frame of $ \Mcal $. It is an orthonormal basis of $ \R^n $, and we denote it by $ \{ \bm{e}_1,\cdots,\bm{e}_n \} $. For any $ r \in \{ 1,\dots,n \} $, $ \bm{e}_r $ is a linear combination of $ \ga^{(1)},\dots,\ga^{(r)} $, where the coefficient functions are smooth, and in particular the coefficient function of $ \ga^{(r)} $ is nonvanishing. In addition, the frame satisfies the Frenet-Serret equations
	\[
	\begin{cases}
		\bm{e}_1^{'}(t) = \kappa_1(t) \bm{e}_2(t), \\
		\bm{e}_r^{'}(t) = -\kappa_{r-1}(t) \bm{e}_{r-1}(t) + \kappa_r(t) \bm{e}_{r+1}(t), \quad r=2,\dots,n-1, \\
		\bm{e}_n^{'}(t) = -\kappa_{n-1}(t) \bm{e}_{n-1}(t),
	\end{cases}
	\]
	where $ \kappa_1,\dots,\kappa_{n-1} $ are smooth positive functions.
	
	Now if $ \bm{j} \in S_1(R,\lambda) $, then for some $ t \in [-2,2] $, we have
	\[
		|\bm{j} \cdot \bm{e}_r(t)|
		\ll \sum_{i=1}^{r} |\bm{j} \cdot \ga^{(i)}(t)|
		\ll \sum_{i=1}^{r} q^{-1} \lambda^i
		\ll q^{-1} \lambda^r,
		\quad r=1,\dots,n.
	\]
	Therefore, there exists a constant $ C_0 > 0 $ depending only on $ \Mcal $ and $ n $ such that
	\begin{multline*}
		S_1(R,\lambda) \subset S_2(R,\lambda) \\
		\stackrel{\text{def}}{=} \bigcup_{t \in [-2,2]} \{ \bm{j} \in \Z^n : |q \bm{j}| < 2R, \, |\bm{j} \cdot \bm{e}_r(t)| < C_0 q^{-1} \lambda^r, \, r=1,\dots,n \}.
	\end{multline*}
	
	Set $ \bm{g} = \bm{e}_n $. By the Frenet-Serret equations, $ \bm{g}^{(r)} $ is a linear combination of $ \{ \bm{e}_{n-r},\cdots,\bm{e}_n \} $, where the coefficient functions are smooth, and in particular the coefficient function of $ \bm{e}_{n-r} $ is nonvanishing. Thus, $ \bm{g}(t), \bm{g}^{(1)}(t),\dots,\bm{g}^{(n-1)}(t) $ are linearly independent for any $ t \in [-2,2] $.
	
	For each fixed $ \bm{\xi} \in \R^n $, assume that
	\[
		\bm{\xi} = \sum_{r=0}^{n-1} \rho_r(t) \bm{g}^{(r)}(t), \quad t \in [-2,2].
	\]
	Since $ \bm{g}^{(r)} \cdot \bm{e}_{n-r} $ is smooth and nonvanishing, and $ \bm{g}^{(r)} \cdot \bm{e}_i = 0 $ when $ i < n-r $, we have
	\[
	\begin{bmatrix}
		\bm{\xi} \cdot \bm{e}_1 \\
		\vdots                  \\
		\bm{\xi} \cdot \bm{e}_n 
	\end{bmatrix}
	=
	\begin{bmatrix}
		0					  & \dots   & \bm{g}^{(n-1)} \cdot \bm{e}_1     \\
		\vdots				  & \iddots & \vdots                            \\
		\bm{g} \cdot \bm{e}_n & \dots   & \bm{g}^{(n-1)} \cdot \bm{e}_n
	\end{bmatrix}
	\begin{bmatrix}
		\rho_{0}   \\
		\vdots     \\
		\rho_{n-1}
	\end{bmatrix},
	\]
	where the matrix is lower anti-triangular and nonsingular.
	
	Now if $ \bm{\xi} \in S_2(R,\lambda) $, using the fact that the inverse of a lower anti-triangular matrix is an upper anti-triangular matrix, for some $ t \in [-2,2] $ we have
	\[
		|\rho_r(t)|
		\ll \sum_{i=1}^{n-r} |\bm{\xi} \cdot \bm{e}_i(t)|
		\ll \sum_{i=1}^{n-r} q^{-1} \lambda^i
		\ll q^{-1} \lambda^{n-r},
		\quad r=0,\dots,n-1.
	\]
	Therefore, enlarging the constant $ C_0 $ if necessary, we get
	\begin{equation}
		S(R,\lambda) \subset S_2(R,\lambda) \subset \Z^n \cap G(R,\lambda), \label{subset:S}
	\end{equation}
	where we define
	\[
		G(R,\lambda) = \bigcup_{t \in [-2,2]} \left\lbrace \sum_{r=0}^{n-1} \rho_r \bm{g}^{(r)}(t) : |\rho_r| < C_0 \, l_r(R,\lambda), \, r=0,\dots,n-1 \right\rbrace,
	\]
	where
	\[
		l_r(R,\lambda) = \max \{ q^{-1} \min{\{ R,\lambda^{n-r} \}}, 1 \}.
	\]
	
	\subsection{Volume bound}
	When $ \lambda > q $, we use the trivial bound:
	\begin{equation}
		\# S(R,\lambda) \ll (q^{-1} R)^n. \label{ineq:S_trivial}
	\end{equation}
	Now we consider the case $ \lambda \leq q $, in which $ l_{n-1}(R,\lambda) = 1 $. Using (\ref{subset:S}) we may reduce the problem of estimating the cardinality of $ S(R,\lambda) $ to that of estimating the volume of $ G(R,\lambda) $.
	
	For brevity, sometimes we omit the argument list $ (R,\lambda) $. For any set $ Z \subset \R^n $, we denote its $ \sigma $-neighbourhood by $ Z^\sigma $ and its $ n $-dimensional Lebesgue measure by $ |Z| $.
	
	Define
	\[
		G_0 = \bigcup_{t \in [-2,2]} \left\lbrace \sum_{r=0}^{n-2} \rho_r \bm{g}^{(r)}(t) : |\rho_r| < C_0 \, l_r, \, r=0,\dots,n-2 \right\rbrace.
	\]
	Since $ l_{n-1} = 1 $, there exists a constant $ \sigma_0 > 0 $ depending only on $ \Mcal $ and $ n $, such that $ G \subset G_0^{\sigma_0} $. Set $ \sigma = \sigma_0 + 1/2 $. Then we have $ G^\frac12 \subset G_0^{\sigma} $. Therefore,
	\begin{equation}
		\# S \asymp |S^\frac12| \leq |G^\frac12| \leq |G_0^\sigma|. \label{ineq:S}
	\end{equation}
	Now it suffices to estimate $ |G_0^{\sigma}| $.
	
	Define a set $ \Omega = [-2,2] \times C_0 [-l_0, l_0] \times \cdots \times C_0 [-l_{n-2}, l_{n-2}] $ and a map
	\[
		\bm{\phi} : \Omega \rightarrow \R^n, \quad (t,s_0,...,s_{n-2}) \mapsto \sum_{r=0}^{n-2} s_r \bm{g}^{(r)}(t),
	\]
	where $ C_0 [-l_r, l_r] = [- C_0 l_r, C_0 l_r] $ for $ r=0,\dots,n-2 $. Then $ G_0^{\sigma} \subset \bm{\phi}(\Omega)^{\sigma} $. It is easy to see that
	\[
		\sum_{r=0}^{n-1} |\bm{g}^{(r)}(t)| \ll 1, \quad t \in [-2,2].
	\]
	Since $ l_r \leq q^{-1} R $ for all $ r = 0,\dots,n-2 $, we have
	\begin{equation}
		\left|  \frac{\partial \bm{\phi}}{\partial t}(\bm{u}) \right| \ll q^{-1} R, \,
		\left|  \frac{\partial \bm{\phi}}{\partial s_r}(\bm{u}) \right| \ll 1, \quad
		\bm{u} \in \Omega, \, r = 0,\dots,n-2. \label{ineq:phi}
	\end{equation}
	
	Suppose that $ \bm{x} \in \bm{\phi}(\Omega)^{\sigma} $. Then there exists some $ \bm{u} = (t,s_0,...,s_{n-2}) \in \Omega $, such that $ |\bm{x} - \bm{\phi}(\bm{u})| < \sigma $. By (\ref{ineq:phi}), for $ \bm{u}' = (t',s_0',...,s_{n-2}') \in \Omega $ satisfying $ |t' - t| \ll q R^{-1}, \, |s_r' - s_r| \ll 1, \, r = 0,\dots,n-2 $, we have $ |\bm{\phi}(\bm{u}') - \bm{\phi}(\bm{u})| < \sigma $ and thus $ |\bm{x} - \bm{\phi}(\bm{u}')| < 2\sigma $. Therefore, noting that $ l_r \geq 1 $ for all $ r = 0,\dots,n-2 $, we have
	\begin{equation}
		\int_{\Omega} \chi_{B(\bm{0},2\sigma)}(\bm{x} - \bm{\phi}(\bm{u})) \, \dd \bm{u} \gg
		q R^{-1} \, \chi_{\bm{\phi}(\Omega)^{\sigma}}(\bm{x}), \quad \bm{x} \in \R^n, \label{ineq:convolution}
	\end{equation}
	where $ \chi_{B(\bm{0},2\sigma)} $ is the characteristic function of the ball $ B(\bm{0},2\sigma) $ in $ \R^n $ centered at the origin with radius $ 2\sigma $. Integrating both sides of (\ref{ineq:convolution}) and using the Fubini-Tonelli theorem, we get
	\[
		|B(\bm{0},2\sigma)| \cdot |\Omega| \gg q R^{-1} |\bm{\phi}(\Omega)^{\sigma}|.
	\]
	Recall that both $ \sigma $ and $ C_0 $ depend only on $ \Mcal $ and $ n $, so
	\[
		|G_0^{\sigma}| \leq |\bm{\phi}(\Omega)^{\sigma}| \ll q^{-1} R \prod_{r=0}^{n-2} l_r.
	\]
	Then by (\ref{ineq:S}) we finally get
	\begin{equation}
		\# S(R,\lambda) \ll q^{-1} R \prod_{r=0}^{n-2} l_r(R,\lambda), \quad \lambda \leq q. \label{ineq:S final}
	\end{equation}
	
	\section{Final calculation} \label{section:final calculation}
	Define $ V(R,\lambda) $ to be the right-hand side of (\ref{ineq:S final}) and $ W(R,\lambda) = V(R,\lambda) \cdot \lambda^{-1} $. That is,
	\begin{align*}
		W(R,\lambda)
		&= \lambda^{-1} q^{-1} R \prod_{r=0}^{n-2} l_r(R,\lambda) \\
		&= \lambda^{-1} q^{-1} R \prod_{r=0}^{n-2} \max \{ q^{-1} \min{\{ R,\lambda^{n-r} \}}, 1 \}.
	\end{align*}
	Let $ x = x(q,\delta) $ satisfy $ q^x = C \delta^{-1} q^{\varepsilon} $ (here we assume $ q > 1 $ without loss of generality), which is equivalent to $ \delta = C q^{\varepsilon - x} $. By  (\ref{ineq:Sigma}), (\ref{ineq:S_trivial}) and (\ref{ineq:S final}), we have
	\begin{align}
		\varSigma_\Mcal(q,\delta)
		&\ll \delta^{n-1} q
		\sum_{\substack{q \leq R < q^{1+x} \\ q^{-1} R \; \text{dyadic}}}
		\sum_{\substack{q < \lambda < C R \\ \lambda \; \text{dyadic}}}
		(q^{-1} R)^n \lambda^{-1}  \notag \\
		&\quad + \delta^{n-1} q
		\sum_{\substack{q \leq R < q^{1+x} \\ q^{-1} R \; \text{dyadic}}}
		\sum_{\substack{C^{-1} R^\frac{1}{n} < \lambda \leq q \\ \lambda \; \text{dyadic}}}
		W(R,\lambda). \label{ineq:Sigma<0+1}
	\end{align}
	The first term on the right-hand side of (\ref{ineq:Sigma<0+1}) is easy to estimate:
	\begin{align}
		\varSigma_0(q,\delta)
		&\stackrel{\text{def}}{=}
		\delta^{n-1} q
		\sum_{\substack{q \leq R < q^{1+x} \\ q^{-1} R \; \text{dyadic}}}
		\sum_{\substack{q < \lambda < C R \\ \lambda \; \text{dyadic}}}
		(q^{-1} R)^n \lambda^{-1} \notag \\
		&\ll \delta^{n-1} q^{nx}. \label{ineq:Sigma_0}
	\end{align}
	As for the second term of the right hand side of (\ref{ineq:Sigma<0+1}), we split the inner sum into $ n $ ranges:
	\begin{align}
		\varSigma_1(q,\delta) 
		&\stackrel{\text{def}}{=}
		\delta^{n-1} q
		\sum_{\substack{q \leq R < q^{1+x} \\ q^{-1} R \; \text{dyadic}}}
		\sum_{\substack{C^{-1} R^\frac{1}{n} < \lambda \leq q \\ \lambda \; \text{dyadic}}}
		W(R,\lambda) \notag \\
		&\leq \delta^{n-1} q
		\sum_{\substack{q \leq R < q^{1+x} \\ q^{-1} R \; \text{dyadic}}}
		\sum_{k=1}^{n} \varLambda_k(R), \label{ineq:Sigma_1}
	\end{align}
	where
	\[
	\varLambda_k(R) =
	\begin{cases}
		\sum_{\substack{R^\frac12 \leq \lambda \leq q \\ \lambda \; \text{dyadic}}} W(R,\lambda),
		& \text{for} \; k = 1, \\
		\sum_{\substack{R^\frac{1}{k+1} \leq \lambda < R^\frac{1}{k} \\ \lambda \; \text{dyadic}}} W(R,\lambda),
		& \text{for} \; k = 2,3,\dots,n-1, \\
		\sum_{\substack{C^{-1} R^\frac{1}{n} < \lambda < R^\frac{1}{n} \\ \lambda \; \text{dyadic}}} W(R,\lambda),
		& \text{for} \; k = n.
	\end{cases}
	\]
	Now we need to estimate $ \varLambda_k(R) $ for $ k = 1,\dots,n $.
	
	If $ R^\frac12 \leq \lambda \leq q $, then
	\[
		(l_0,\dots,l_{n-2}) = (q^{-1} R,\dots,q^{-1} R), \; V(R,\lambda) = V(R,R^\frac12) = q^{-n} R^n.
	\]
	Thus,
	\begin{equation}
		\varLambda_1(R) = V(R,R^\frac12)
		\sum_{\substack{R^\frac12 \leq \lambda \leq q \\ \lambda \; \text{dyadic}}} \lambda^{-1}
		\ll V(R,R^\frac12) \cdot R^{-\frac12} = W(R,R^\frac12). \label{ineq:Lambda_1}
	\end{equation}
	
	If $ C^{-1} R^\frac{1}{n} < \lambda < R^\frac{1}{n} $, note that for fixed $ R $, the function $ V(R,\lambda) $ is non-decreasing in $ \lambda $, then we have
	\begin{equation}
		\varLambda_n(R) \leq V(R,R^\frac{1}{n})
		\sum_{\substack{C^{-1} R^\frac{1}{n} < \lambda < R^\frac{1}{n} \\ \lambda \; \text{dyadic}}} \lambda^{-1}
		\ll V(R,R^\frac{1}{n}) \cdot R^{-\frac{1}{n}} = W(R,R^\frac{1}{n}). \label{ineq:Lambda_n}
	\end{equation}
	
	Now we consider the more complicated case $ R^\frac{1}{k+1} \leq \lambda < R^\frac{1}{k} $, $ k \in \{ 2,\dots,n-1 \} $. In this case, we have
	\[
		(l_0,\dots,l_{n-2}) = (q^{-1} R, \dots, q^{-1} R, \max \{ q^{-1} \lambda^k, 1 \}, \dots, \max \{ q^{-1} \lambda^2, 1 \}).
	\]
	
	If $ q^{-1} \lambda^k > 1 $, which is equivalent to $ \lambda > q^\frac{1}{k} $, then for fixed $ R $, the function $ V(R,\lambda) $ grows at least as fast as $ \lambda^2 $. That is to say, for $ \lambda > q^\frac{1}{k} $ and fixed $ R $, the function $ V(R,\lambda) \cdot \lambda^{-2} = W(R,\lambda) \cdot \lambda^{-1} $ is non-decreasing in $ \lambda $. 
	
	On the other hand, if $ q^{-1} \lambda^k \leq 1 $, which is equivalent to $ \lambda \leq q^\frac{1}{k} $, then for fixed $ R $, the function $ V(R,\lambda) $ is constant. By the above analysis, we have
	\begin{align}
		\varLambda_k(R)
		&\leq \sum_{\substack{R^\frac{1}{k+1} \leq \lambda \leq q^\frac{1}{k} \\ \lambda \; \text{dyadic}}}
		V(R,\lambda) \cdot \lambda^{-1} +
		\sum_{\substack{q^\frac{1}{k} < \lambda < R^\frac{1}{k} \\ \lambda \; \text{dyadic}}}
		W(R,\lambda) \cdot \lambda^{-1} \cdot \lambda \notag \\
		&\leq V(R,R^\frac{1}{k+1})
		\sum_{\substack{R^\frac{1}{k+1} \leq \lambda \leq q^\frac{1}{k} \\ \lambda \; \text{dyadic}}}
		\lambda^{-1} +
		W(R,R^\frac{1}{k}) \cdot R^{-\frac{1}{k}}
		\sum_{\substack{q^\frac{1}{k} < \lambda < R^\frac{1}{k} \\ \lambda \; \text{dyadic}}}
		\lambda \notag \\
		&\ll V(R,R^\frac{1}{k+1}) \cdot R^{-\frac{1}{k+1}} +
		W(R,R^\frac{1}{k}) \cdot R^{-\frac{1}{k}} \cdot R^\frac{1}{k} \notag \\
		&= W(R,R^\frac{1}{k+1}) + W(R,R^\frac{1}{k}). \label{ineq:Lambda_k}
	\end{align}
	Combining (\ref{ineq:Lambda_1}), (\ref{ineq:Lambda_n}) and (\ref{ineq:Lambda_k}), we get
	\begin{equation}
		\sum_{k=1}^{n} \varLambda_k(R) \ll \sum_{k=2}^{n} W(R, R^\frac{1}{k}). \label{ineq:Lambda sum 1}
	\end{equation}
	
	According to the definition, for $ k \in \{ 2,\dots,n \} $ and $ R \leq q^{1 + \frac{1}{n-1}} $, we have
	\[
		W(R, R^\frac{1}{k}) = R^{-\frac{1}{k}} (q^{-1} R)^{n-k+2}
		\prod_{i=2}^{k-1} \max \{ q^{-1} R^\frac{i}{k}, 1 \}
		= R^{-\frac{1}{k}} (q^{-1} R)^{n-k+2}.
	\]
	Substituting (\ref{ineq:Lambda sum 1}) into (\ref{ineq:Sigma_1}), we obtain, for $ x \leq 1/(n-1) $,
	\begin{align}
		\varSigma_1(q,\delta) 
		&\ll \delta^{n-1} q
		\sum_{k=2}^{n} \sum_{\substack{q \leq R < q^{1+x} \\ q^{-1} R \; \text{dyadic}}}
		R^{-\frac{1}{k}} (q^{-1} R)^{n-k+2} \notag \\
		&\ll \delta^{n-1} q \sum_{k=2}^{n} q^{(n-k+2-\frac1k)x - \frac1k}. \label{ineq:Sigma_1 final}
	\end{align}
	Combining (\ref{ineq:Sigma<0+1}), (\ref{ineq:Sigma_0}) and (\ref{ineq:Sigma_1 final}), for $ x \leq 1/(n-1) $, we have
	\begin{equation}
		\varSigma_\Mcal(q,\delta) \ll \delta^{n-1} q \sum_{k=1}^{n} q^{(n-k+2-\frac1k)x - \frac1k}. \label{ineq:Sigma final}
	\end{equation}
	
	Let $ x_k $ be the unique solution to $ (n-k+2-\frac1k)x - \frac1k = 0 $, that is,
	\[
		x_k = \frac{1}{k(n-k+2)-1}, \quad k = 1,2,\dots,n.
	\]
	A simple optimization gives
	\[
		x_0 \stackrel{\text{def}}{=} \min_{1 \leq k \leq n} x_k = \frac{1-\Theta(n)}{n-1},
	\]
	where $ \Theta(n) $ is as defined in Theorem~\ref{thm:hickman&srivastava}.
	
	Therefore, if $ x < x_0 $, then $ \delta^{n-1} q \gg 1  $ and by (\ref{ineq:Nfra}) and (\ref{ineq:Sigma final}) we have
	\[
		\Nfra_\Mcal(q,\delta) \ll \delta^{n-1} q + \OO_\varepsilon(1) \ll_\varepsilon \delta^{n-1} q.
	\]
	If $ x \geq x_0 $, which is equivalent to $ \delta \leq C q^{\varepsilon - x_0} $, then by the increasing property of $ \Nfra_\Mcal(q,\delta) $ with respect to $ \delta $, we have
	\begin{align*}
		\Nfra_\Mcal(q,\delta)
		&\leq \Nfra_\Mcal(q,C q^{\varepsilon - x_0}) \\
		&\ll (q^{\varepsilon - x_0})^{n-1} q \, (1 + \sum_{k=1}^{n} q^{(n-k+2-\frac1k)x_0 - \frac1k}) + \OO_\varepsilon(1) \\
		&\ll_\varepsilon q^{\Theta(n) + (n-1)\varepsilon}. 
	\end{align*}
	Combining the preceding estimates and choosing $ \omega = \omega^+ $, $ \eta = \eta^+ $ and $ \varepsilon = \nu / (n-1) $, we get
	\[
		\Nfra_{\omega^+,\eta^+,\Mcal}(q,\delta) \ll_\nu \delta^{n-1} q + q^{\Theta(n) + \nu}
		\quad \text{for } q \geq 1 \text{ and } \delta \in (0,\frac12).
	\]
	Then by (\ref{ineq:N}) we have proved the first part of Theorem~\ref{thm:reduced}.
	
	Assume that $ x \leq x_0 - \varepsilon $, or equivalently, that $ \delta \geq C q^{(\Theta(n)-1) / (n-1) + 2\varepsilon} $. Then $ \delta^{n-1} q \gg q^\frac12  $ and by (\ref{ineq:Sigma final}) we have $ \varSigma_\Mcal(q,\delta) \ll \delta^{n-1} q^{1-\varepsilon} $. So by (\ref{ineq:Nfra}) we have
	\[
		|\Nfra_\Mcal(q,\delta) - c_0(q) \delta^{n-1} q| \leq C_1 \delta^{n-1} q^{1-\varepsilon} + C_\varepsilon,
	\]
	where $ C_1 $ depends only on $ \Mcal $ and $ n $, and $ C_\varepsilon $ depends only on $ \varepsilon $, $ \Mcal $ and $ n $. Then for $ q \gg_\varepsilon 1 $, we have
	\[
		C_1 q^{-\varepsilon} < \frac13 c_0(q) \text{ and } C_\varepsilon < \frac13 c_0(q) \delta^{n-1} q
	\]
	and thus
	\[
		\Nfra_\Mcal(q,\delta) > \frac13 c_0(q) \delta^{n-1} q.
	\]
	Taking $ \omega = \omega^- $, $ \eta = \eta^- $ and $ \varepsilon = \nu / 2 $ in the preceding argument, we obtain
	\[
		\Nfra_{\omega^-,\eta^-,\Mcal}(q,\delta) \gg \delta^{n-1} q
		\quad \text{for } q \gg_\nu 1 \text{ and } \frac12 > \delta \gg q^{\frac{\Theta(n)-1}{n-1} + \nu}.
	\]
	Then by (\ref{ineq:N}) we have proved the second part of Theorem~\ref{thm:reduced}.
	
	\subsection*{Acknowledgments}
	Helpful discussions and support from Prof. Jingwei Guo are gratefully acknowledged. This work was partially supported by the NSFC (Grant Nos. 12571110 and 12341102). The author thanks Rajula Srivastava for pointing out an error in a previous version of this article.
	
	\subsection*{AI usage}
	Yuanbao and ChatGPT were used to assist with the English-language presentation of this paper.

\end{document}